\documentclass[10pt]{amsart}

\usepackage{amsmath,amssymb,stmaryrd,dsfont}
\usepackage[pdftex]{graphicx}
\usepackage[all]{xy}
\usepackage{etaremune}
\usepackage{pdfsync}

\newcommand{\var}{\varepsilon}

\newcommand{\e}{\varepsilon}

\newcommand{\be}{\begin{equation}}
\newcommand{\beq}{\begin{equation}}
\newcommand{\ee}{\end{equation}}
\newcommand{\bac}{\begin{array}{c}}
\newcommand{\ea}{\end{array}}
\begin{document}
\newtheorem{theorem}{Theorem}[section]
\newtheorem{lemma}[theorem]{Lemma}
\newtheorem{observation}[theorem]{Observation}
\newtheorem{definition}[theorem]{Definition}
\newtheorem{example}[theorem]{Example}
\newtheorem{corollary}[theorem]{Corollary}
\newtheorem{claim}[theorem]{Claim}
\newtheorem{assumption}{Assumption}
\newtheorem{property}{Property}

\numberwithin{equation}{section}

\title[Selection of  classical limit]{On the selection of the classical limit for potentials with BV derivatives}

\author[Athanassoulis]{Agissilaos  ATHANASSOULIS}
\address[A.  Athanassoulis]{Department of Applied Mathematics, University of Crete,  P.O. Box 2208,
71409 Heraklion, GREECE}
\email{athanassoulis@tem.uoc.gr}
\author[Paul]{Thierry Paul}
\address[T. Paul]{CNRS and CMLS, \'Ecole polytechnique, 91128 Palaiseau cedex, FRANCE}
\email{thierry.paul@math.polytechnique.fr}

%\date{}

\today

%\xxivtime

\begin{abstract} 
We consider the classical limit of the quantum evolution, with some rough potential, of wave packets concentrated near singular trajectories of the underlying dynamics. We prove that under appropriate conditions, even in the case of BV vector fields, the correct classical limit can be selected.

\end{abstract}

\maketitle
\tableofcontents
%
%
%\begin{itemize}
% \item
%\end{itemize}
%
%\newpage

\section{Introduction}

The fundamental equation of quantum mechanics is formulated either as the Schr\"{o}dinger equation for a wavefunction,
\begin{equation}
\label{eq1}
\begin{array}{c}
i\varepsilon \frac{\partial }{\partial t}u^\varepsilon(t) 
= \left( -\frac{ \varepsilon^2}{2}\Delta +V\left( x \right) \right)
u^\varepsilon(t) , \\
u^\varepsilon(t=0)=u^\varepsilon_0(x),
\end{array}
\end{equation}
or more generally as the Heisenberg-von Neumann equation for a density matrix   $D^\e(t)$,
\begin{equation}
\label{eq11}
\begin{array}{c}
i\varepsilon \frac{\partial }{\partial t}D^\varepsilon(t) = \left[ -\frac{ \varepsilon^2}{2}\Delta +V,
D^\varepsilon(t)\right] , \\
D^\varepsilon(t=0)=D^\varepsilon_0.
\end{array}
\end{equation}
The state of the quantum system is described by the operator $D^\e$ (understood to be the orthonormal projector $|u^\e \rangle \langle u^\e |$ when working with the Schr\"odinger equation), evolving in time under the aforementioned equations. Indeed, under very general conditions on the potential $V$, known as Kato's conditions,  this is a problem well posed for $u^\e_0 \in L^2$ or $D^\e_0$ being a Hilbert-Schmidt operator \cite{kat0,kat1/2,kat1}. It is also well known that the positivity and trace of $D^\varepsilon_0$ is preserved in time. The parameter $\e$ is called Planck's constant, and when it is very small one usually expects the system to behave like a classical one.

An equivalent way to write equation (\ref{eq11}) is in terms of the Wigner transform,
\begin{equation}\label{density}
 W^\varepsilon(x,k,t)=W^\e[D^\e(t)](x,k)=\int\limits_{y \in \mathbb{R}^n} { e^{-2\pi i yk} K^\varepsilon
 (x+\varepsilon \frac{y}2,x-\varepsilon \frac{y}2,t) dy }
\end{equation}
where for each time $t$ $K^\varepsilon(x,y,t)$ is the integral kernel of the operator $D^\e(t)$; in other words, for $f \in \mathcal{S}$
\[
D^\e f (x) = \e^{-n} \int{e^{2\pi i \frac{x-y}{\e}k} W^\e(\frac{x+y}2,k) f(y)dy dk}.
\]
A compact way to say that is that the Wigner transform of an operator is $\e^n$ times its Weyl symbol. Thus the operator corresponding to the Wigner function $W^\e_0$ is $D^\e_0=\e^{-n} Op_{Weyl} \left({ W^\e_0 }\right)$.

 In the case of a {\em pure state} this yields
\begin{equation}
\label{eqyhnju}
 W^\e(x,k,t)=W^\e[u^\e(t)](x,k)=\int\limits_{y \in \mathbb{R}^n} { e^{-2\pi i yk}u^\varepsilon (x+\varepsilon \frac{y}2,t) 
 \bar{u}^\varepsilon (x-\varepsilon \frac{y}2,t) dy }.
\end{equation}
Equation (\ref{eq11}) implies \cite{LP}
\begin{equation}
\label{eq1m2aol}
\begin{array}{c}
\partial_t {{W}^\varepsilon} +  2\pi k \cdot \partial_x {{W}^\varepsilon} %+ \\ { } \\

+i\int{e^{-2\pi i S y} \frac{V(x+\frac{\e}{2}y)-V(x+\frac{\e}{2}y)}{\var} dy \,\, W^\e(x,k-S,t)dS}=0,

\\ { } \\
{W}^\varepsilon(t=0)={W}^\varepsilon_0=W^\e[D^\e_0].
\end{array}
\end{equation}
The propagator for equation (\ref{eq1m2aol}) is constructed from the one of equation (\ref{eq11}) and as long as Kato's conditions apply and $D^\e_0$ is a Hilbert-Schmidt operator (equivalently $W_0^\e \in L^2$) there exists a unique solution for problem (\ref{eq1m2aol}), and it is the Wigner transform of the solution of equation (\ref{eq11}); see theorem \ref{thrmWignL2} below, and  \cite{Mark}.

We will often use the shorthand $T^V_\e$ for the operator involving the potential (see Definitions and Notations in section \ref{secDef} below).

The Wigner transform (WT) respects the structure of the problem, allowing the Wigner equation (\ref{eq1m2aol}) to inherit properties (and in particular conservation laws) from (\ref{eq11}), see e.g. lemma \ref{thrmWignL2}.
Moreover, under appropriate conditions, it allows for a very natural and compact description of the semiclassical limit, $\e \shortrightarrow 0$. Indeed the WT has a physically meaningful limit as $\var$ tends to zero, while in general the wavefunction $u^\varepsilon$ itself (or the operator $D^\e$) does not. The limit (in the weak-$*$ topology of an appropriate algebra of test functions), called the Wigner measure,
\[
W^\varepsilon \rightharpoonup W^0, 
\]
satisfies the Liouville equation of classical  mechanics
\begin{equation}
\label{eqlolLIOUV}
\begin{array}{c}
\partial_t W^0 + 2\pi k \partial_x W^0 - \frac{1}{2\pi} \partial_x V \cdot \partial_k W^0=0, \\ { } \\
W^0(t=0) =\mathop{lim}\limits_{\varepsilon \shortrightarrow 0} W^\varepsilon_0
\end{array}
\end{equation}
if the potential is regular enough \cite{LP,ger}.

For potentials with low regularity (less than $C^{1,1}$)  two different problems naturally arise;

{\bf Problem 1:} showing that the Wigner measure is a weak solution of (\ref{eqlolLIOUV}); 

\noindent and, if possible, 

{\bf Problem 2:} constructing a selection principle that identifies the correct weak solution, since problem (\ref{eqlolLIOUV}) can be ill-posed in that case. 

\vskip 0.2cm
Problem 1 is solved in \cite{LP} for $V \in C^1$ (here we will work for quite less regular, if somewhat specific, potentials). As is highlighted already for $V \in C^1 \setminus C^{1,1}$, there can be several weak solutions when the initial datum is a singular measure (i.e. supported on a set of measure zero, the most natural case being a single point). That is what we mean in the sequel when we say the quantum initial data $W^\e_0$ is {\em concentrating} --    concentrating to a point-supported measure $W^0_0 = \mathop{lim}\limits_{\var \shortrightarrow 0} W^\e_0$. Indeed, the Liouville equation can be made well-posed for much worse potentials than $C^1$, but typically not for concentrating initial data \cite{ambrosio,ldp,gmmp,Mark}.

More recently, Problem 1 has been  successfully handled in \cite{AFFGP,FLP} under BV condition of the vector field generated by the classical Hamiltonian, and for \textbf{non-concentrating} (in the sense defined above) and {\bf small} (in operator norm) initial  data. In that situation the solution of \eqref{eq1m2aol} tends weakly to the push forward of the limit of the initial datum by the so-called DiPerna-Lions-Bouchut-Ambrosio flow (defined a.e.). In \cite{AP} a new technique was formulated for handling semiclassical limits with concentrating initial data. Parts of it were used to work out problems with low-regularity potentials (also in strong topology),  in \cite{apbad}; see also \cite{tp} for a  short review. In a nutshell, the technique of \cite{apbad} applies to more general initial data than the results of \cite{AFFGP,FLP} (including in particular data {\bf concentrating to a point-supported measure}), but less general potentials (roughly speaking, of $C^{1,a}$ type).

In the present note  we compute explicitly the Wigner measures for certain problems with potentials that are not in $C^1$, but have measure valued second derivative (therefore in a sense having the typical flavour of BV vector fields), and initial data which concentrate on the points of singularity. We illustrate the idea behind this computation and discuss how  this idea can be generalized to other problems with similar regularity, and, roughly speaking, isolated repulsive singularities. The main difference from \cite{apbad} is that here we manage to work with other topologies, more appropriate for the problem, and this is crucial to achieving stronger results. 

The main result of this paper is theorem \ref{thrm1}; let us also mention the simple observation of Lemma \ref{lmcortrnrm}.
%We have a technique that can select the correct semiclassical limit in some cases, and at the moment it's not clear how to characterize in an intuitive, economical way the class of problems to which it applies. We  illustrate the proof, point out some exotic constructions that can be achieved with it, and outline how to check its applicability to a given problem (and what kind of problems can one reasonably hope that it applies to). A more precise and conclusive characterization of its applicability will have to be given in a forthcoming publication.

\section{Main result }

%The main result of our paper will be given in the form of an example. We will discuss extensions in section \ref{ext}.
% Although the potential $V$ in the following result is fixed, let us mention that the same result can be obtained for general potential having the same type of singularity by replacing $\psi$ below by any compactly supported smooth function non vanishing at the origin (?). 

Since the regularity condition on the potential does not insure unicity of the flow, we will compare the solution of the quantum problem \eqref{eq1m2aol} to the behaviour of the solution of the Liouville equation \textbf{with the same, $\var$-dependent, initial condition}. We will show that the limit of the solution of the  problem \eqref{eq1m2aol}  is the same as the limit of the solution of the Liouville equation with the quantum (i.e. $\var$-dependent) initial datum. In other words, the extra information needed for the selection of the weak solution that correctly captures the semiclassical limit is fully contained in the way in which the initial datum concentrates to a point-supported measure.

\begin{theorem} \label{thrm1} Let $\psi: [0,\infty) \shortrightarrow [0,1]$ be a $C^\infty$ monotone cutoff function, $\psi(s)=1$ for all $s<\frac{1}2$, $\psi(s)=0$ for all $s>1$. Now, for $x \in \mathbb{R}^2$, some $\theta \in [0,1)$ set
\[
V(x)=\left({ 1-|x_1|^{1+\theta} }\right)\,\, \psi(|x_1|)\psi(|x_2|).
\]
Moreover let
\[
F_0^\e = 
 \delta_x^{-2} \delta_k^{-2} w\left({\frac{x-(0,-X)}{\delta_x},\frac{k-(0,K)}{\delta_k}}\right) 
\]
for some $w \in H^2 \cap L^\infty \cap W^{1,1}$, $w(x,k) \geqslant 0$, $\int{wdxdk}=1$, $supp \,\, w \, \subseteq \{ x^2+k^2 \leqslant 1 \}$, $X>1,K>0$ and 
  $\delta_k =\delta_x^2=(-log(\e^{\frac{1}4}))^{\frac{1}{\theta-2}}$ (see claim \ref{claimncc} for other possible scalings).
\vskip 0.2cm
Let ${W}^\varepsilon$ be the solution of
\begin{equation}
\label{eq1hsdvaol}
\begin{array}{c}
\partial_t {{W}^\varepsilon}(x,k,t) +  2\pi k \cdot \partial_x {{W}^\varepsilon}(x,k,t) +  \\

+i\int{e^{-2\pi i S y} \frac{V(x+\frac{\e}{2}y)-V(x+\frac{\e}{2}y)}{\var} dy \,\, W^\e(x,k-S,t)dS}=0,

\\ { } \\
{W}^\varepsilon(t=0)=F_0^\e.
\end{array}
\end{equation}
(Recall that $W^\e$ is the WT of the solution of equation (\ref{eq11}) with initial data $D^\e_0=\e^{-n} Op_{Weyl} \left({ W^\e_0 }\right)$).
Denote also by $\rho^\e$ the solution of
\beq \label{eqziouzid}  
\bac
\partial_t \rho^\e + 2\pi k\cdot \partial_x \rho^\e -\frac{1}{2\pi} \partial_x V(x) \cdot \partial_k \rho^\e=0, \\ { } \\

\rho^\e(t=0)=F^\e_0.
\ea
\ee
(With this regularity of $V$ problem (\ref{eqziouzid}) has a unique solution in $L^p$ \cite{ambrosio}).

\vskip 0.4cm
Then, for each $t \in [0,T]$
\beq \label{eqsmigg}
\mathop{lim}\limits_{\e \shortrightarrow 0} \langle W^\e(t) - \rho^\e(t), \phi \rangle =0 \hskip 1cm \forall \phi \in \mathcal{A} \cap L^2 .
\ee
(See Definition \ref{defap} below for the algebra $\mathcal{A}$).

In other words  the semiclassical limit of (\ref{eq1hsdvaol}) behaves exactly in the same way as the concentration limit for the Liouville equation (\ref{eqziouzid}).
\end{theorem}

\vskip 0.3cm
It can be noted that the only kind of initial data allowed in theorem \ref{thrm1} are wavepackets concentrating on a given point. 
One reason the result is phrased the way it is, is that approximation (\ref{eqsmigg}) can hold even when there is no unique semiclassical limit. If, for example we restrict the values of $\e \in \{\frac{1}m\}$ and we substitute the initial datum $F_0^\e=F_0^{\frac{1}m}$  by $F_0^{\frac{1}m}(x_1-d^{\frac{1}m},x_2,k_1,k_2)$ for $d^{\frac{1}m}=C(-1)^m (log \, m)^{-\frac{1}2}$, then it is easy to check that the result applies, but $W^{\frac{1}m}(t)$ --    and $\rho^{\frac{1}m}(t)$ with it --    has two accumulation points, one scattered to the right of $\{x_1=0\}$ and the other to the left (see also next remark). The finding here is that, whatever the interaction with the singularity, it is {\em the same} for the quantum problem (\ref{eq1hsdvaol}) and the classical problem (\ref{eqziouzid}). In other words, all the information needed to determine the interaction is contained in the initial datum. In that light, even this limited pool of initial data contains enough different possibilities to explore.

\vskip 0.3cm
We can use the theorem to compute in more detail the semiclassical limit:
\begin{corollary}
\be\label{split}
W^0(t)=c_+\delta{(\mathcal{X}^+(t),\mathcal{P}^+(t))}+c_-\delta{(\mathcal{X}^-(t),\mathcal{P}^-(t))},
\ee
with $c_\pm = \int\limits_{\pm x_1 >0}{w(x,k)dxdk}$
and  $(\mathcal{X}^\pm(t),\mathcal{P}^\pm(t))$ are the two diverging curves obtained by $(\mathcal{X}^\pm(t),\mathcal{P}^\pm(t))=\mathop{lim}\limits_{\eta \shortrightarrow 0} (\mathcal{X}^\pm_\eta(t),\mathcal{P}^\pm_\eta(t))$,
\[
\bac
\dot{\mathcal{X}}_\eta^\pm(t)=2\pi \mathcal{P}_\eta^\pm(t), \,\,\,\,\,\, \dot{\mathcal{P}}_\eta^\pm(t)=-\frac{1}{2\pi} \nabla V(\mathcal{X}_\eta^\pm(t)), \\
\mathcal{X}^\pm_\eta(t=0)=(\pm \eta,-X), \,\,\,\,\, \mathcal{P}_\eta^\pm(t=0)=(0,K).
\ea
\]
\end{corollary}

\begin{figure}[htb!]
\centering
\includegraphics[width=110mm,height=75mm]{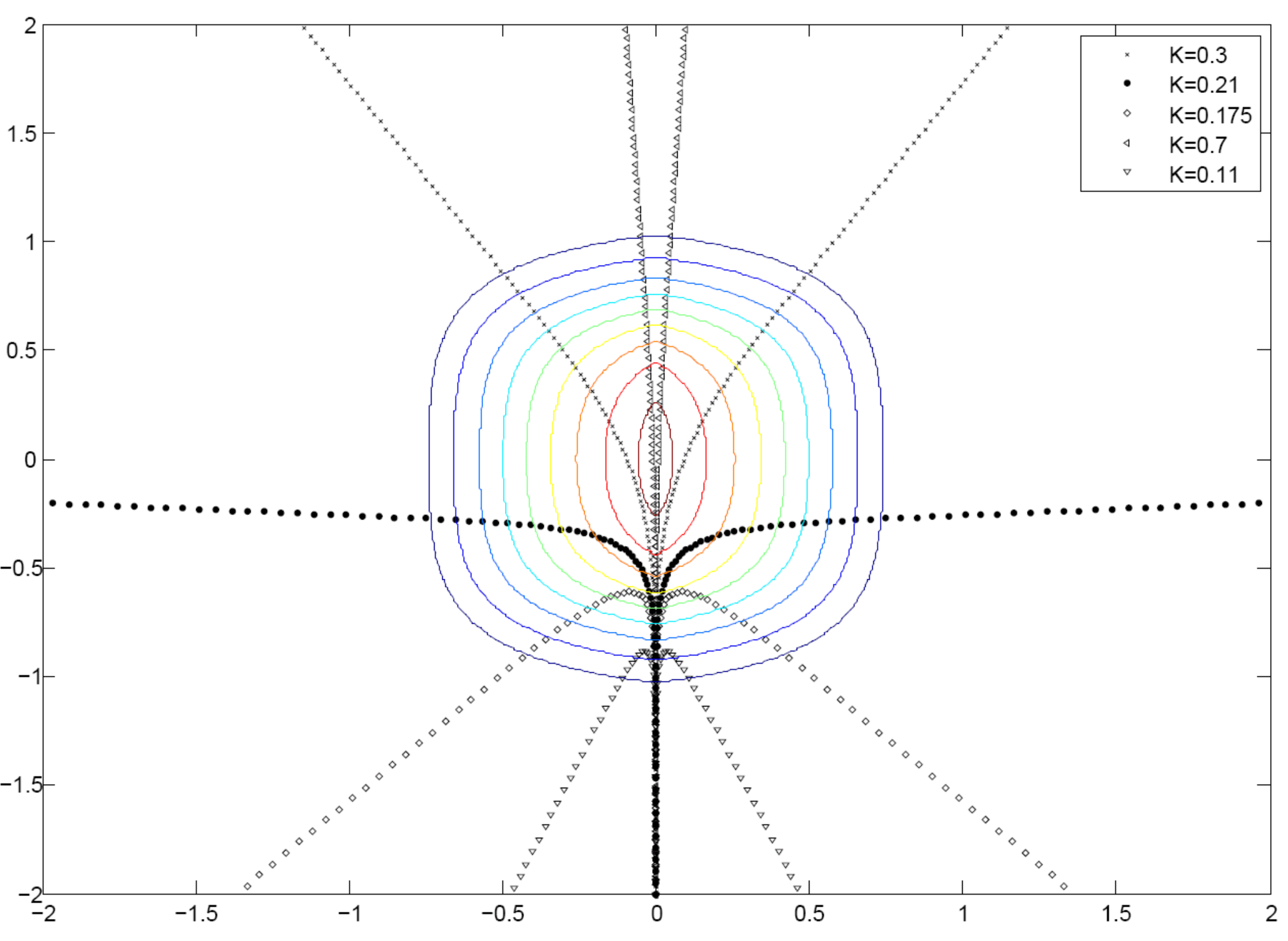}
\caption{ 
Pairs of trajectories $(\mathcal{X}^\pm,\mathcal{P}^\pm)$ corresponding to various values of the parameter $K$. The projection to the $(\mathcal{X}_1,\mathcal{X}_2)$ plane, as well as the contours of the potential $V(\mathcal{X}_1,\mathcal{X}_2)$ are shown.
} 
\label{Fig2201}
\end{figure}

This can in fact be used to construct quite exotic examples; in figure \ref{Fig2201} one sees various trajectories  (their projection on the $\mathcal{X}$ plane, to be precise) for various values of the parameter $K$. (The value of the parameter $X$ is not really interesting, since it merely determines the length of free motion before interaction with the potential starts). Observe moreover, that by an appropriate superposition (discrete or continuous) of initial data as in theorem \ref{thrm1}, one can construct an example where the outgoing waves are (discretely or continuously) distributed among almost any angle.

This is a highly exotic example, and it may well turn out that its physical relevance is limited. However, it does point out a couple of points that could prove rather fruitful. First of all, the tools of semiclassical analysis can be extended and applied even to phenomena that are qualitatively very different from the smooth regime that they originate in. Such splitting of a particle and other possible exotic examples can motivate technical (and possibly numerical) experiments and refinements on existing methods and tools.

Moreover, this example points quite naturally to a more general (and, one expects, physical) situation: there is a ``scattering  process'' happening ``on the singularity''. That is, scattering with respect to local fast variables a neighbourhood of the singularity. In this particular case (of very slow concentration, since $\delta_x,\delta_k$ are logarithmic in the semiclassical parameter) the scattering operator is ``trivially'' given by the concentration limit of the respective classical problem. That is, it is an operator fully determined and constructed by classical dynamics. In cases of faster concentration (e.g. pure states) one expects that an auxiliary {\em quantum problem} will need to be solved in a neighbourhood of the singularity to determine the ``splitting''. In other words the irregular potential makes the corresponding flow to behave in $O(1)$ time-scales as a regular one would in {\em long times} (diverging trajectories), and this comparison can yield some very insightful finds \cite{tp2}.

\vskip 0.3cm
It should also be noted that
this particular example rests on the exact alignment of the direction of propagation to the line of singularity (even a small rotation  can, potentially, destroy it). As can be easily checked, if the potential is substituted by $V(x)=\left({ 1-|x|^{1+\theta} }\right)\,\, \psi(|x_1|)\psi(|x_2|)$ theorem \ref{thrm1} still holds, and then is robust to any perturbation of the initial data. Still, the only {\em nontrivial} case is when the particle passes {\em exactly} over the singularity, and the possibility of non-standard effects in that case. 

At this point we must refer once again to \cite{AFFGP}; a stable result, which covers ``almost all initial data'' --    in which case the a.e. unique solution of the Liouville equation describes the semiclassical limit, and there is no room for exotic effects. In other words these interactions may be ``atypical'' in some sense, but if one tries to study them, they {\em have} to select initial data so that they fully interact.

\vskip 0.2cm
Let us finally remark that an outline of how to check whether this technique applies to different problems is given in section \ref{ext}.

%
%
%
%\vskip 1cm
%Finally, it should be noted that
%this result generalizes our previous work \cite{apbad} in two ways; firstly it applies to potentials with measure valued second derivatives (in particular which may not be in $C^1$); moreover by significantly weakening the ``non-concentration condition'', here stated and shown to hold in claim \ref{claimncc} (That is, even with $\theta >0 \Rightarrow V \in C^1$, \cite{apbad} would not apply here) and finally by strating with an initial datum which is not pined up on the singularity.
%

\section{Definitions and Notations}%\label{def}
\label{secDef}

The Fourier transform is defined as
\begin{equation}
\widehat{f}(k)=\mathcal{F}_{x \shortrightarrow k} \left[ {f(x)} \right]
=\int\limits_{x\in \mathbb{R}^n} {e^{-2\pi i k x} f(x)dx}. 
\end{equation}

 For compactness, we will use the following notations:
\be
\mathcal{F}_2 W(x,K)=\int{e^{-2\pi i k K}W(x,k)dk}.
\ee
\vskip 0.1cm
\be \label{eqwop}
\begin{array}{c}
T^V_\var W = i\int{e^{-2\pi i S y} \frac{V(x+\frac{\e}{2}y)-V(x+\frac{\e}{2}y)}{\var} dy \,\, W(x,k-S)dS}=
\\ { } \\

=\mathcal{F}^{-1}_{K \shortrightarrow k} \left[{ 
\frac{V(x+\frac{\e}{2}K)-V(x-\frac{\e}{2}K)}{\var} \mathcal{F}_{2} W(x,K)
 }\right] %= \\ { } \\
= \mathcal{F}^{-1}_{K \shortrightarrow k} \left[{ 
\int\limits^{\frac{K}{2}}_{s=-\frac{K}{2}} { \nabla V(x+s)ds } \,\, \mathcal{F}_{2} W(x,K)
 }\right],
\end{array}
\ee
\vskip 0.1cm
\be
\begin{array}{c}
T^V_0 W=-\frac{1}{2\pi} \partial_x V \cdot \partial_k W= \mathcal{F}^{-1}_{K \shortrightarrow k} \left[{ 
\nabla V(x) \cdot K \,\, \mathcal{F}_{2} W(x,K)
 }\right].
\end{array}
\ee

\begin{definition} \label{defap} The space $\mathcal{A}$ is defined as the completion of the smooth functions of compact support $C^\infty_c(\mathbb{R}^{2n})$ under the norm
\[
||f||_{\mathcal{A}} =\int\limits_{K}{ \mathop{sup}\limits_{x} |\mathcal{F}_2 f(x,K)|dK}.
\]
It follows that
\[
||f||_{\mathcal{A}'} =\mathop{sup}\limits_{K} \int\limits_{x}{  |\mathcal{F}_2 f(x,K)|dx}.
\]
Some basic properties that we will use are
\begin{eqnarray} \nonumber
||f||_{L^\infty} & \leqslant & ||f||_{\mathcal{A}} \\
||f||_{\mathcal{A}'} & \leqslant & ||f||_{L^1} \leqslant ||f||_{\mathcal{M}}
\nonumber
\end{eqnarray}
\end{definition}

\begin{definition} \label{defmes} Every finite signed measure can be decomposed to positive and negative part,
$
\mu = \mu^+ -\mu^-,
$
for some finite non-negative measures $\mu^+,\mu^-$. We will denote the total variation of a signed measure
\[
||\mu ||_{\mathcal{M}} = \int\limits_{\mathbb{R}^n} d\mu^+ + \int\limits_{\mathbb{R}^n} d\mu^- .
\]
\end{definition}

\begin{definition} \label{defPHI}
We will denote by $\Phi=\Phi^\e$ the smoothing operator
\[
\Phi^\e : f(x) \mapsto \left({\frac{2}{\e}}\right)^{\frac{n}{2}} \int{ e^{-2\pi \frac{(x-x')^2}{\e}} f(x')dx'}.
\]
When there is no danger of confusion, we will write
\[
\widetilde{f}:=\Phi^\e f .
\]
\end{definition}

\begin{definition} The Sobolev space $W^{m,p}$ is defined as the completion   of the smooth functions of compact support $C^\infty_c(\mathbb{R}^{n})$ under the norm
\[
||f||_{W^{m,p}} =\sum\limits_{|A|\leqslant m} ||\nabla^A f||_{L^p}.
\]
\end{definition}

\section{Proof of the main result}

\subsection{Strategy of the proof}

The intuition behind this proof is quite simple: the quantum initial datum, while concentrating to a point still is an $L^2$ function, hence the a.e. theory for its evolution under the Liouville equation applies. The concentration limit of problem (\ref{eqziouzid}) is therefore a natural candidate for the semiclassical limit. Of course to show a semiclassical approximation, a certain degree of smoothness in the potential and the Wigner function is needed, and in any frontal approach to this problem, such smoothness simply is not there.

A simple idea to try out is the following: can we cut-off the pieces of the Wigner function that approach too closely to the singularity? If we do that, can we meaningfully use an auxiliary function supported just far enough away from the singularity so that it preserves enough smoothness itself --    as well as allowing one to cut-off the potential's singularities? This was pretty much the program we followed before, in \cite{apbad}.

The new element here is that to strengthen that technique to potentials as bad as the ones we treat here, basically a much bigger piece of the initial datum would have to be cut-off, and we need to find a meaningful way in which such a cut-off introduces ``small'' errors. To do that, the positivity of the density matrix comes into play, and some very different considerations are needed.

Claim \ref{claimncc} holds all the compromises that need to be made, and is the conclusion of a lot trial and error. If one assumes it and move on in a first reading, the flow of the rest proof should provide a reasonable motivation for why these computations have to be just so.

\subsection{Proof of theorem \ref{thrm1} }

To facilitate the presentation, let us introduce at this point a number of auxiliary functions (using the notations introduced in the previous section):
\begin{eqnarray}
\label{eq1appsol001}
\partial_t {{W}^\varepsilon} +  2\pi k \cdot \partial_x {{W}^\varepsilon} 
+T^V_\var W^\e=0, \\ 
\nonumber
{W}^\varepsilon(t=0)=F^\varepsilon_0, \\
\nonumber { } \\
\label{eq1appsol002}
\partial_t {W}^\varepsilon_1 +  2\pi k \cdot \partial_x {{W}^\varepsilon_1} 
+T^V_\var W^\e_1=0, \\ 
\nonumber
{W}^\varepsilon_2(t=0)= F^\var_1=\Phi F^\varepsilon_0, \\
\nonumber { } \\
\label{eq1appsol003}
\partial_t {W}^\varepsilon_2 +  2\pi k \cdot \partial_x {{W}^\varepsilon_2} 
+T^V_\var W^\e_2=0, \\ 
\nonumber
{W}^\varepsilon_2(t=0)=F^\var_2=\Phi (1-\psi(\frac{x}{R'})) F_0^\var, \\
\nonumber { } \\
\label{eq1appsol004}
\partial_t {W}^\varepsilon_3 +  2\pi k \cdot \partial_x {{W}^\varepsilon_3} 
+T^V_\var W^\e_3=0, \\ 
\nonumber
{W}^\varepsilon_3(t=0)=F^\var_3=(1-\psi(\frac{x}{R'})) F_0^\var, \\
\nonumber { } \\
\label{eq1appsol005}
\partial_t \rho^\varepsilon_1 +  2\pi k \cdot \partial_x \rho^\varepsilon_1 
+T^V_0 \rho^\e_1=0, \\ 
\nonumber
\rho^\varepsilon_1(t=0)=F^\var_3, \\
\nonumber { } \\
\label{eq1appsol005}
\partial_t \rho^\varepsilon +  2\pi k \cdot \partial_x \rho^\varepsilon 
+T^V_0 \rho^\e=0, \\ 
\nonumber
\rho^\varepsilon(t=0)=F^\var_0.
\end{eqnarray}
The function $\psi$ is the same as in the statement of theorem \ref{thrm1}, and $R'$ will be set below.

Obviously, we are going to use these functions as stepping stones, passing from one to the other with the appropriate topology each time. Because this topology cannot be always the same, the end result is formulated as is, in weak sense: for all $\phi \in L^2 \cap \mathcal{A}$
\[
\bac
|\langle W^\var(t)-\rho^\var(t), \phi \rangle| \leqslant |\langle W^\var(t)-W_1^\var(t), \phi \rangle|+
|\langle W_1^\var(t)-W_2^\var(t), \phi \rangle|+ \\
+|\langle W_2^\var(t)-W_3^\var(t), \phi \rangle|+
|\langle W_3^\var(t)-\rho_1^\var(t), \phi \rangle|+
|\langle \rho_1^\var(t)-\rho^\var(t), \phi \rangle|
\ea
\]
The point is to collect the various constraints that would come from each of these building-block problems and satisfy them at the same time. That is essentially done in the following

\begin{claim} \label{claimncc} With appropriate calibration of $\delta_x,\delta_k,R'$, we have 
\[\bac ||F_3^\e-F_0^\e||_{L^1}=o(1), \\ { } \\ 

||\rho_1^\e||_{L^\infty([0,T],H^2)}=o(\e^{-\frac{1}{2}}),\ea\]
and
\[
||F_2^\e-F_3^\e||_{L^2}=o(1).
\]
\end{claim}

\noindent {\bf Proof of the claim:}
The point is to make sure that $\rho_1^\e(t)$ doesn't pass through a neighbourhood of the set where the second and third derivatives of the potential are singular (or too large in any case); in this case the strip $\{ |x|<R \}$, $R=(-log(\e^{\frac{1}4}))^{\frac{1}{\theta-2}}$, for $t \in [0,T]$. In that case, using lemma \ref{thrmRegLiouv} it follows readily that indeed 
\beq
\label{eqtybzann}
||\rho_1^\e||_{L^\infty([0,T],H^2)}=O(\e^{-\frac{1}{4}} ||F_3^\e||_{H^2}).
\ee
(by recalling that in fact only the values of the potential's derivatives along the part of phase space that the solutions passes through matter).

To ensure that $\rho_1^\e(t)$ stays away from $\{ |x|<R \}$, it turns out to be sufficient to exclude a somewhat larger strip $\{ |x|<R' \}$. To see why, imagine firstly that we are in free space, $V=0$. The largest possible momentum in the $x_1$ direction is $\delta_k$; in time $T$ this can only cover a distance of $T\delta_k$; therefore $R'=R+T\delta_k$ should do it.

\begin{figure}[htb!]
\centering
\includegraphics[width=45mm,height=48mm]{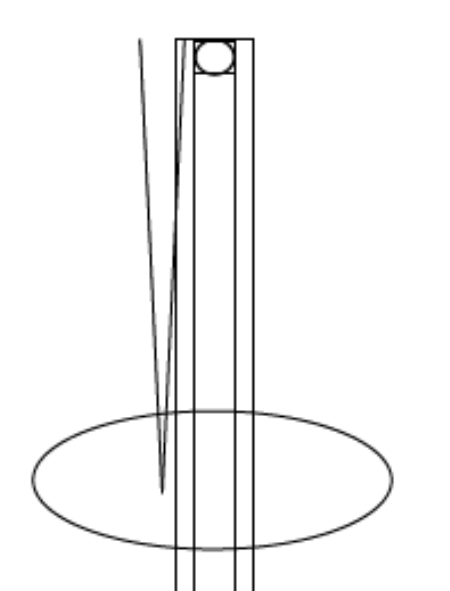}
\caption{ Any trajectory leaving the support of $F_3^\varepsilon$ has an initial velocity in a small cone around $k=(0,K)$. If we were in free space, there wouldn't be enough time for sufficient movement in $x_1$ to reach $\{|x_1|<R\}$ for $t \in [0,T]$. Therefore the solution is the sum of two components, supported on either side of $\{|x_1|<R\}$. The time-scale $T$ enters in the constants used, but the construction can be made for any fixed $T$.} 
\label{Fig222}
\end{figure}

Now we have to take into account the addition of the potential. The repulsive nature of the singularity (preserved by the cutoff $\psi(|x_1|)$) makes sure that any trajectory moving towards it will not be accelarted --    but in fact slowed down; i.e. any movement towards the singularity will be {\em less} than $T\delta_k$. The presence of $\psi(|x_2|)$ serves only to slow somewhat movement in the $x_2$ direction (or even turn it back); but makes no difference whatsoever in the $x_1$ direction. 

The timescale $T$ is interesting when trajectories starting in the support of $F_0^\e$ reach and leave the support of the potential; in any case $T=O(1)$.

\vskip 0.2cm
Moreover, for the first part of the claim,
\[
\bac
||F_3^\e-F_0^\e||_{L^1} \leqslant ||(1-\Phi^\e)\delta_x^{-2} \delta_k^{-2} w\left({\frac{x}{\delta_x},\frac{k}{\delta_k}}\right)||_{L^1}+\\
+ ||\psi(\frac{x}{R'})\delta_x^{-2} \delta_k^{-2} w\left({\frac{x}{\delta_x},\frac{k}{\delta_k}}\right)||_{L^1} \leqslant \\ { } \\

\leqslant \sqrt{\e} ||\delta_x^{-2} \delta_k^{-2} w\left({\frac{x}{\delta_x},\frac{k}{\delta_k}}\right)||_{W^{1,1}} + ||w||_{L^\infty} ||\psi(\frac{x \delta_x }{R'})||_{L^1}= \\ { } \\

=O(\frac{ \sqrt{\e}}{\delta_x}+\frac{ \sqrt{\e}}{\delta_k}  + \frac{R'}{\delta_x})
\ea
\]

 \vskip 0.2cm 
 Finally, for the third part,
 \[
||F_2^\e-F_3^\e||_{L^2}=||(I-\Phi^\e)F_3^\e||_{L^2} \leqslant \sqrt{\e} ||F_3^\e||_{H^1} = O( \frac{\sqrt{\e}}{R} (\delta_x^{-1}+\delta_k^{-1})\delta_x^{-1}\delta_k^{-1}  ) 
\]
 
\vskip 0.5cm 
So collecting all the constraints, we have
\[
\bac
\sqrt{\e} \ll \delta_x, \delta_k; \,\, R'\ll \delta_x  \,\,\,\, \Rightarrow \,\,\,\, 
||F_3^\e-F_0^\e||_{L^1} =o(1),  \\ { } \\

R'= (-log(\e^{\frac{1}4}))^{\frac{1}{\theta-2}}  + C \delta_k \,\,\,\, \Rightarrow \,\,\,\, ||\rho_1^\e||_{L^\infty([0,T],H^2)}=O(\e^{-\frac{1}{4}} ||F_3^\e||_{H^2}), \\ { } \\

(\delta_x^{-2} + \delta_k^{-2})\delta_x^{-1}\delta_k^{-1} \ll \e^{-\frac{1}{4}} \,\,\,\, \Rightarrow \,\,\,\, ||F_3^\e||_{H^2} = o(\e^{-\frac{1}{4}}), \\ { } \\

\sqrt{\e} \ll R' \delta_x^{2}\delta_k, \,\, \sqrt{\e} \ll R' \delta_x\delta_k^{2} \,\,\,\, \Rightarrow \,\,\,\, ||F_2^\e-F_3^\e||_{L^2}=o(1).
\ea
\]
A concrete scaling that makes all these constraints covalid is (recall that $R=(-log(\e^{\frac{1}4}))^{\frac{1}{\theta-2}}$)
\[
\delta_k=\delta_x^2=R.
\]
The proof of the claim is complete.

\vskip 0.2cm
\noindent {\bf Remark: } The claim apparently depends on the geometry of the problem; but really all we used was that the flow is repulsive away from the singularity. That is, that any trajectory would not go towards the line $\{ x_1=0 \}$ faster than it would on free space. Thus it makes no difference if e.g.  a different potential --    for which the same property holds --    is used. Indeed if the potential was $V(x)=\left({ 1-|x|^{1+\theta} }\right)\,\, \psi(|x_1|)\psi(|x_2|)$ nothing needs to be changed.

\vskip 0.5cm

Given this calibration, the proof proceeds in a very predictable fashion. Indeed, using the unitary propagation of the Wigner equation --  theorem \ref{thrmWignL2} --  one observes that
\[
||F_0^\var-F_1^\var||_{L^2} \leqslant \sqrt{\var} ||F_0^\var||_{H^1} =O(\sqrt{\var} \delta_x\delta_k(\delta_x+\delta_k) )= o(1) \,\, \Rightarrow \,\, ||W_0^\var(t)-W_1^\var(t)||_{L^2}=o(1) \,\, \forall t.
\]
\vskip 0.2cm

The estimate between $W_1^\var$ and $W_2^\var$ is one of the essential parts, and the technical innovation here. Denote by $U^\e(t)$ the propagator of the Wigner equation; we observe that 
\[
\bac
W_1^\e - W^\e_2=U^\e(t)(F_1^\e-F_2^\e)=\\

=U^\e(t)\left({ \frac{2}{\e} }\right)^2 
\int{  e^{-2\pi \frac{(x-x')^2+(k-k')^2}{\e}} \, \psi(\frac{x'}{R}) \delta_x^{-2} \delta_k^{-2} w\left({\frac{x'-(0,-2)}{\delta_x},\frac{k'-(0,1)}{\delta_k}}\right) \, dx'dk'},
\ea
\]
the point being that $F_1^\e-F_2^\e$ is itself a density matrix --  convolution of a positive measure with a coherent state; see lemma \ref{dktnkmmtl}. Moreover, it is a {\em small} density matrix, and this is preserved in time; see lemma \ref{eqvoss}.
Indeed it follows by the conservation of trace (and lemma \ref{lmcortrnrm} for density matrices) that 
 \[|| W_1^\e(t) - W^\e_2(t)||_{\mathcal{A}'}=|| F^\e_1 - F^\e_2||_{\mathcal{A}'}\leqslant|| F^\e_1 - F^\e_2||_{L^1} =O(\frac{R'}{\delta_x})
 \]
which was already calibrated to be $o(1)$ in the proof of claim \ref{claimncc}. (Indeed following that calibration we have $\frac{R'}{\delta_x}=O(\sqrt{\delta_k}+\delta_x)$). Observe that $F^\e_1 - F^\e_2$ is {\em not} $o(1)$ in $L^2$ sense, hence the introduction of the $L^1$-like norm is necessary.

\vskip 0.2cm
Moreover, by virtue of theorem \ref{thrmWignL2} and claim \ref{claimncc}
\[
||W_2^\e(t)-W^\e_3(t)||_{L^2}=||F_2^\e-F_3^\e||_{L^2} =o(1).
\]

\vskip 0.2cm
The other non-trivial step is the quantum-classical dynamics comparison between $W^\e_3$ and $\rho_1^\e$.
Recall that $U^\e(t)$ the propagator of the Wigner equation, and denote by $E(t)$ the propagator of the Liouville equation.
Set  $h^\e(t)=W^\e_3(t)-\rho_1^\e(t)$; then of course $h^\e(t) =  \int\limits_{\tau=0}^t{U^\e(t) (T^{{V}}_\e-T^{{V}}_0) \rho_1^\e(\tau)  d\tau}$.

Using the Fourier transform in the $k$ variable we have
\[
\bac
||\int\limits_{\tau=0}^t{U(t) (T^{{V}}_\e-T^{{V}}_0) \rho_2^\e(\tau)  d\tau}||_{L^2} \leqslant \\ { } \\

\leqslant T \mathop{sup}\limits_{t \in [0,T]} || ( \frac{ {V}(x+\frac{\e K}{2})-{V}(x-\frac{\e K}{2}) }{\e} -\partial_x {V}(x) \cdot K) \mathcal{F}_2 \rho_2^\e(t) ||_{L^2} \leqslant \\ { } \\

\leqslant T \mathop{sup}\limits_{t \in [0,T]}  || \chi_{[0,\e^{-\frac{1}{2}}]}(|K|) ( \frac{ {V}(x+\frac{\e K}{2})-{V}(x-\frac{\e K}{2}) }{\e} -\partial_x {V}(x) \cdot K) \mathcal{F}_2 \rho_2^\e(t) ||_{L^2} + \\ { } \\

+T \mathop{sup}\limits_{t \in [0,T]}|| \chi_{(\e^{-\frac{1}{2}},+\infty)}(|K|) ( \frac{ {V}(x+\frac{\e K}{2})-{V}(x-\frac{\e K}{2}) }{\e} -\partial_x {V}(x) \cdot K) \mathcal{F}_2 \rho_2^\e(t) ||_{L^2}  \leqslant \\ { } \\

\leqslant  \e || \chi_{[0,\e^{-\frac{1}{2}}]}(|K|)    \frac{ |K|^2}{2} \mathop{sup}\limits_{
\begin{scriptsize}
\begin{array}{c}
\tau \in (-1,1) \\
|A|=2
\end{array}
\end{scriptsize}
} | \partial_x^A {V}( x+ \tau\frac{\e K}{2}) | \,\, \mathcal{F}_2 \rho_2^\e(t)||_{L^2} + O(\e^\frac{1}{2}) ||\, |K|^2 \mathcal{F}_2 \rho_2^\e(t)||_{L^2}.
\ea
\]

For the first term we used a straightforward Taylor expansion, while for the complementary case we used the computation
\[
\bac
|( \frac{ {V}(x+\frac{\e K}{2})-{V}(x-\frac{\e K}{2}) }{\e} -\partial_x {V}(x) \cdot K)| \leqslant \\ { } \\

\leqslant |\frac{1}{\e} \int\limits_{-\frac{\e|K|}{2}}^{\frac{\e|K|}{2}} (\frac{\e|K|}{2} -t) \sum\limits_{j,m} z_j z_m \partial_{x_j x_m} V(x+t z) dt| \leqslant C |K| \sum\limits_{j,m} |\int\limits_{-\frac{\e|K|}{2}}^{\frac{\e|K|}{2}}   \partial_{x_j x_m} V(x+t z) dt| \leqslant \\ { } \\

\leqslant C |K| \sum\limits_{A=2} ||\partial_x^A V||_{\mathcal{M}}.
\ea
\]

The key to proceed is to make use of the well prepared initial datum $F_3^\e$: denote
\[
S^\e=\{ x | \exists \,\, t \in [0,T], k \in \mathbb{R}^2 \,\, : \,\, \rho_1^\e(x,k,t) > 0 \}
\]
and
\[
B^\e= \mathop{\bigcup}\limits_{x_* \in S^\e} \{ x| |x-x_*|< \sqrt{\e} \}.
\]
Then
\[
\chi_{[0,\e^{-\frac{1}{2}}]}(|K|)    \mathop{sup}\limits_{
\begin{scriptsize}
\begin{array}{c}
\tau \in (-1,1) \\
|A|=2
\end{array}
\end{scriptsize}
} | \partial_x^A {V}( x+ \tau\frac{\e K}{2}) | \leqslant \mathop{sup}\limits_{
\begin{scriptsize}
\begin{array}{c}
x \in B^\e \\
|A|=2
\end{array}
\end{scriptsize}
} | \partial_x^A {V}( x) |
\]
and therefore
\[
\bac
\e  \chi_{[0,\e^{-\frac{1}{2}}]}(|K|)    \mathop{sup}\limits_{
\begin{scriptsize}
\begin{array}{c}
\tau \in (-1,1) \\
|A|=2
\end{array}
\end{scriptsize}
} | \partial_x^A {V}( x+ \tau\frac{\e K}{2}) | \,\,\,\, || \,\,|K|^2 \mathcal{F}_2 \rho_2^\e(t)||_{L^2} \leqslant \\ { } \\

\leqslant  \e \mathop{sup}\limits_{
\begin{scriptsize}
\begin{array}{c}
x \in B^\e \\
|A|=2
\end{array}
\end{scriptsize}
} | \partial_x^A {V}( x) | \,\, \mathop{sup}\limits_{t \in [0,T]} ||\rho_2^\e(t)||_{H^2}
\ea
\]
finally yielding
\[
\bac
||\int\limits_{\tau=0}^t{U(t) (T^{{V}}_\e-T^{{V}}_0) \rho_2^\e(\tau)  d\tau}||_{L^2} = \\ { } \\

=O(\sqrt{\e} ||\rho_2^\e(t)||_{H^2} (1+ \sqrt{\e} \mathop{sup}\limits_{
\begin{scriptsize}
\begin{array}{c}
x \in B^\e \\
|A|=2
\end{array}
\end{scriptsize}
} | \partial_x^A {V}( x) |) ) = O(\e^{\frac{1}4} ||F_3^\e||_{H^2} (1+R^{\theta-1}))=o(1).
\ea
\]

This is possible precisely because our approximate initial datum $F_3^\e$ stays away from a neighbourhood of the singular set $\{ |x_1|=0 \}$ when pushed forward in time by the a.e. Ambrosio-Lions-Di Perna flow $E(t)$, as was checked in claim \ref{claimncc}.
(For a potential of the form $|x|$ the second derivatives are zero almost everywhere, and this creates the possibility of taking advantage of a very special structure. That's why we included the case $\theta>0$ in the statement of the theorem, to show that in principle this technique work for a variety of localized repulsive singularities).

\section{Auxiliary results}

\begin{lemma}[$2^{nd}$ order derivatives equations for the Liouville equation] \label{thrmRegLiouv}
Consider the Cauchy problem for the Liouville equation with potential $V(x)$ on $\mathbb{R}^n$,
\begin{equation}
\label{eqlolL1I2O3U34V}
\begin{array}{c}
\partial_t f + 2\pi k \cdot \partial_x f - \frac{1}{2\pi} \partial_x V(x) \cdot \partial_k f=0, \\ { } \\
f(t=0) =f_0.
\end{array}
\end{equation}
There are constants $C_1,C_2>0$ depending only on $n$  such that
\[
||f(t)||_{H^2} \leqslant C_1 e^{ tC_2 \mathop{sup}\limits_{|a|\leqslant 3} ||\partial_x^a V(x)||_{L^\infty} }||f_0||_{H^2}.
\]
\end{lemma}

\begin{lemma}[Density matrices]\label{dktnkmmtl} Let $\mu$ be a probability measure on $\mathbb{R}^{2n}$. Then
\[
W_0^\var(x,k)=\left({\frac{2}\var}\right)^n \int{ e^{-2\pi\frac{(x-x')^2+(k-k')^2}{\var}} d\mu(x',k')}
\]
is the Wigner function of a density matrix, i.e. the corresponding operator $D^\var$ is a positive trace-class operator with $tr(D^\var)=1$.
\end{lemma}

\noindent {\bf Proof:} Though the proof is obvious by using T\"oplitz quantization, let us give a direct proof.

We know that $tr(D^\var)=\int{W^\var dxdk}=1$. Let us now look at positivity. To that end, observe that
the integral kernel
\[
K^\var(x,y)=\left({\frac{2}\var}\right)^{\frac{n}2} \int\limits_{x_0,k_0}{ e^{2\pi i\frac{k_0}{\var}(x-y)}  e^{-\frac{\pi}{\var}[(x-x_0)^2+(y-x_0)^2]} d\mu(x_0,k_0) }
\]
is the kernel of a positive operator. Indeed:
\[
\begin{array}{c}
\int{K^\var(x,y)u(x)\overline{u}(y)dxdy} = \left({\frac{2}\var}\right)^{\frac{n}2} \int\limits_{x_0,k_0}{ |\langle e^{-2\pi i \frac{k_0}\var x+\frac{\pi}\var (x-x_0)^2} ,u \rangle|^2 d\mu(x_0,k_0) } \geqslant 0.
\end{array}
\]

The proof is complete by observing that the Wigner function corresponding to the kernel $\rho^\var$ is
\[
\int{e^{-2\pi i ky} K^\var(x+\frac{\var y}2,x-\frac{\var y}2)dy}=W_0^\var.
\]

\begin{lemma} \label{lmcortrnrm} For any trace class operator $D^\e$, with corresponding Wigner function $W^\e$,
\[
||W^\e||_{\mathcal{A}'}\leqslant ||D^\e||_{tr}=tr(|D^\var|).
\]
Moreover, if $D^\e \geqslant 0$
\[
||W^\e||_{\mathcal{A}'}=||D^\e||_{tr}=tr(D^\var).
\]
\end{lemma}

\proof Let $D^\e$ be a non-negative trace-class operator. Then, it admits a SVD expansion over orthonormal projectors
\[
D^\e=\sum\limits_m \lambda_m |u_m\rangle \langle u_m|,
\]
where of course $||\lambda_m||_{l^1}=||D^\e||_{tr}=tr(D^\e)$, i.e. $\lambda_m \geqslant 0$, and $\langle u_m, u_l \rangle =\delta_{m,l}$. It follows that
\[
W^\e=\sum\limits_m \lambda_m W^\e[u_m]=\sum\limits_m \lambda_m \int{e^{-2\pi i k y} u_m(x+\frac{\e y}2)\overline{u}_m(x-\frac{\e y}2)dy}
\]
and, by straightforward substitution,
\[
\bac
||W^\e||_{\mathcal{A}'}= \mathop{sup}\limits_{K} \int\limits_{x}{  |\sum\limits_m \lambda_m u_m(x+\frac{\e K}2)\overline{u}_m(x-\frac{\e K}2)|dx}  \leqslant \\ { } \\ 

\leqslant\sum\limits_m \lambda_m \mathop{sup}\limits_{K} \int\limits_{x}{    |u_m(x+\frac{\e K}2)\overline{u}_m(x-\frac{\e K}2)|dx}=\sum\limits_m \lambda_m.
\ea
\]
Now observe that
\[
\bac
||W^\e||_{\mathcal{A}'}= \mathop{sup}\limits_{K} \int\limits_{x}{  |\sum\limits_m \lambda_m u_m(x+\frac{\e K}2)\overline{u}_m(x-\frac{\e K}2)|dx}  \geqslant \\ { } \\ 

\geqslant \int\limits_{x}{  \sum\limits_m \lambda_m |u_m(x)|^2 dx} = \sum\limits_m \lambda_m.
\ea
\]
Hence, $D^\e \geqslant 0 \,\, \Rightarrow \,\, ||W^\e||_{\mathcal{A}'}=tr(D^\e)=||D^\e||_{tr}$.

\vskip 0.25cm
To conclude, observe that for $D^\e=D^\e_+ - D^\e_{-}$
\[
|| W^\e||_{\mathcal{A}'} \leqslant ||W^\e[D^\e_+]||_{\mathcal{A}'}  +  ||W^\e[D^\e_-]||_{\mathcal{A}'} =tr(D^\e_+) + tr(D^\e_-)= ||D^\e||_{tr}.
\]

The proof is complete

\vskip 0.2cm

The following theorem is well known, and follows from the correspondence between the density matrix and the Wigner function; see e.g. Theorem 2.1 of \cite{Mark}:  
\begin{theorem}[$L^2$ regularity of the Wigner equation] \label{thrmWignL2} If the Schr\"odinger operator $-\frac{1}2\Delta +V(x)$ is essentially self-adjoint on $L^2(\mathbb{R}^n)$, then the corresponding Wigner equation preserves the $L^2$ norm, i.e. for $W^\varepsilon_0 \in L^2(\mathbb{R}^{2n})$ there is a unique solution of
\begin{equation}
\label{WignerEq}
\begin{array}{c}
\partial_t W^\varepsilon + 2\pi k \cdot \partial_x W^\varepsilon +\frac{2}{\varepsilon} Re \left[{ i\int{e^{2\pi i Sx}\hat{V}(S)W^\varepsilon(x,k-\frac{\varepsilon S}{2})dS} }\right]=0, \\
W^\varepsilon(t=0)=W^\varepsilon_0,
\end{array}
\end{equation}
and $||W^\varepsilon(t)||_{L^2(\mathbb{R}^{2n})}=||W^\varepsilon_0||_{L^2(\mathbb{R}^{2n})}$ $\forall t \in \mathbb{R}$.
\end{theorem}

\vskip 0.2cm
\begin{lemma}[Trace conservation] \label{eqvoss}
Let $W^\e_0$ be a Wigner function corresponding to a trace-class operator. Consider a potential $V$ such that the corresponding Schr\"odinger operator is essentially self-adjoint. If by $W^\e(t)$ we denote its evolution in time under the corresponding Wigner equation, then
\[
tr(W^\e_0)=tr(W^\e(t)) \,\,\,\,\,\,\,\,\, \forall t.
\]
\end{lemma}

\proof One easily checks that $W^\e_0 = \sum\limits_m \lambda_m |u_m\rangle \langle u_m| \Rightarrow W^\e_0 = \sum\limits_m \lambda_m |u_m(t)\rangle \langle u_m(t)|$, where $u_m(t)$ is the evolution in time of $u_m \in L^2$ under the Schr\"odinger equation.

The result follows.

\section{Extensions}\label{ext}

As was mentioned earlier, this approach also applies if the potential is substituted by $V(x)=\left({ 1-|x|^{1+\theta} }\right)\,\, \psi(|x_1|)\psi(|x_2|)$. Other straightforward generalizations come by embedding the problem in higher-dimensional space (i.e. including more transverse dimensions). 

\vskip 0.2cm
Checking whether some version of this approach applies to a problem would start by building the counterpart of claim \ref{claimncc}, and its second part in particular. More specifically:

\begin{itemize}
\item
Set $N=\{ (x,k) | \mathop{sup}\limits_{|A|=3} |\partial_x^A V(x)| > \zeta \}$ ($\zeta$ is a parameter of the order $-log \, \var$ to some power). This is the ``neighbourhood of the singularity'' that we want to stay away from.
\item
Back-propagate it for the appropriate time scale, $M= \bigcup\limits_{t \in [0,T]} \phi_{-t} (N)$. The counterpart to $F_0^3$ should be chosen so that it does not enter $M$. $M$ is not a small set; but we only need that cutting it off makes a small difference to $F_0^\e$, not that itself it is small. $F^\e_0(1-\chi_M(x,k))$ is a starting point for $F^\e_3$; a simpler cutoff might be preferable. Certainly though, if $||F^\e_0(1-\chi_M(x,k))||_{L^1}$ is not small, then this approach does not apply.
\item
The previous step is only one constraint about at least how much we have to cut-off. If we cut-off any less (i.e. exclude a smaller set), we will enter too close to the singularity and the estimates will fail. An other constraint comes from having several $H^2$ bounds to check for our various approximate initial data. This means that if we exclude too small sets (e.g. too thin strips etc), the derivatives of the approximate data will be too large, and the estimates also fail. This is what makes us exclude a strip even when $\theta=0$ in theorem \ref{thrm1}, where the previous step would be fulfilled by taking $N$ to be a line.
\end{itemize}

In other words, there are two different origins for constraints in this construction: on the one hand, we want to avoid the low-regularity region. The more we cut-off the better. On the other hand, we want the cut-offs to introduce small errors, including in various derivative-norms. This is a double sided constraint, as cutting off either a too thick or too thin strip will fail here. Of course the norms we can work in are constrained by the available conservation laws for the equations we work with, i.e. apparently have to be based on $L^2$, and $\mathcal{A}'$ as long as we keep track successfully of positivity. The fact that $\mathcal{A}'$ scales like $L^1$ in concentrating data is crucial in allowing wider strips to be cut-off.

\vskip 0.2cm
The fact that we will have to exclude domains measured in $- (log\, \var)^{-1}$ (to some power) means that in general one can only work with slowly concentrating initial data. Moreover, the {\em repulsive} character of the singularity is used implicitly here, as it helps the set $M$ not be too large, but, when projected on the hyperplane transverse to propagation, basically looks like the set $N$. It is not necessary per se, but a way to meaningfully control how much larger is the set $M$ from $N$ is needed.

\end{document}